\documentclass [a4paper]{article}
\usepackage{amsfonts}
\usepackage{graphicx}
\usepackage{amsmath,amssymb}
\input{epsf.sty}
\usepackage{subfigure}
\usepackage{epsfig}
\usepackage{footmisc}
\usepackage{url}
\usepackage{amsthm}
\usepackage{amscd,enumerate,latexsym}
\usepackage{multirow}
\usepackage{cite}
\usepackage[colorlinks,citecolor=blue,urlcolor=magenta]{hyperref}
\usepackage{breakurl}
\usepackage{booktabs}

\newtheorem{rem}{Remark}
\newtheorem{algorithm}{Algorithm}
\begin{document}
%\renewcommend{theequation}{\thesection.\arabic{equation}}
%\newcommend{\be}{begin{equation}}
%\newcommend{\\ee}{\end{equation}}

\begin{center}
\Large\bf{A flexible and adaptive Simpler GMRES with deflated restarting for shifted linear systems}\\

\quad\\
\normalsize Hong-Xiu Zhong\footnote[1]{School of Science, Jiangnan University, Wuxi, Jiangsu 214122, P.R. China.
E-mail: {\tt zhonghongxiu@126.com}. This author is supported by NSFC (11701225 and 11471122),
the Fundamental Research Funds for the Central Universities (JUSRP11719),
and the Natural Science Foundation of Jiangsu Province (BK20170173).},~
Xian-Ming Gu\footnote[2]{Corresponding author: School of Economic Mathematics/Institute
of Mathematics, Southwestern University of Finance and Economics, Chengdu, Sichuan 611130,
P.R. China. E-mail: {\tt guxianming@live.cn,~guxm@swufe.edu.cn}. This author is supported
by NSFC (1370147, 11501085 and 61402082), the Fundamental Research Funds for the Central
Universities (JBK1809003 and ZYGX2016J138).}

\end{center}

\vspace{0.5cm}

\mbox{\bf Abstract:}~~In this paper, two efficient iterative algorithms based 
on the simpler GMRES method are proposed for solving shifted linear systems. 
To make full use of the shifted structure, the proposed algorithms utilizing 
the deflated restarting strategy and flexible preconditioning can significantly 
reduce the number of matrix-vector products and the elapsed CPU time. Numerical 
experiments are reported to illustrate the performance and effectiveness of the 
proposed algorithms.

\mbox{\bf Keywords: } Shifted linear system, Adaptive Simpler GMRES, Flexible
preconditioning, Deflated restarting.

\mbox{\bf AMS classifications: } 65F15, 65F10, 65Y20.

\section{Introduction}
\setcounter{equation}{0}
\quad\
In this study, we are interested in efficiently simultaneous solutions of the following large shifted linear systems
\begin{equation}
(A+\alpha_jI)x(\alpha_j)=b,~~~~j=1,\cdots,s.
\label{shsystem}
\end{equation}
In general, $A\in\mathbb{C}^{n\times n}$ is non-singular and non-Hermitian, $\alpha_j\in\mathbb{C}$ is the shift such that
$A+\alpha_jI$ is also non-singular, and $\alpha_j$ varies in a wide range, the right-hand side $b\in\mathbb{C}^n$ is fixed.
Usually we take $\alpha_1=0$ as default, otherwise, Eq. (\ref{shsystem}) can be reset after a shift $\alpha_1$. The first
linear system is called the seed system, others are the add systems. Such problem occurs in many scientific and
engineering applications, such as structural dynamics \cite{FerPerSim2000,SimPer2002}, quantum chromodynamics \cite{DarMorWil2004},
web search ranking \cite{LanMey2006}, control theory \cite{DatSad1991,Ahmad2017} and so on. Therefore, there is a strong need for
establishing efficient solutions of Eq. (\ref{shsystem}).

Many traditional methods (such as direct and iterative linear systems solvers)
for the above problem are to solve $(A + \alpha_jI)x(\alpha_j) = b$ for each
$\alpha_j$, this trick can be quite expensive and prohibited when $s$ and $n$ are large. Fortunately, owing to the shift-invariance
property of Krylov subspaces, the Krylov subspace methods can solve Eq. (\ref{shsystem}) simultaneously \cite{Jeger1996}. That is, the Krylov subspace
holds that
\begin{equation*}
\mathcal{K}_m(A,b)=\mathcal{K}_m(A+\alpha_jI,b),~~~~\forall \alpha_j\in\mathbb{C}.
\end{equation*}
Hence, all approximate solutions for (\ref{shsystem}) can be sought in a single
space generated by the matrix $A$ with the vector $b$.

The GMRES algorithm \cite{SaaSch1986} is such a famous Krylov subspace method that
it calculates the basis for $\mathcal{K}_k(A,b)$ by once Arnoldi process with the
initial guess $x_0=0$, hence the shifted system (\ref{shsystem}) can be solved cheaply
if GMRES is performed for it simultaneously \cite{DatSad1991}. However, since the
residuals $r_m(\alpha_j)=b-(A+\alpha_jI)x_m(\alpha_j)$ are not colinear, so that
$\mathcal{K}_m(A,r_m)\neq\mathcal{K}_m(A+\alpha_jI,r_m(\alpha_j))$ with $m$ being
the restarting frequency. As a remedy, Frommer and Gl\"{a}ssner have forced the
residual vectors to be colinear \cite{FroGla1998}, then restarts can again solve
Eq. (\ref{shsystem}) cheaply. There are many variants based on GMRES for solving
shifted linear systems. For instance, Gu, Zhang and Li proposed a variant of the
restarted GMRES augmented with some approximate eigenvectors for the shifted system
(\ref{shsystem}), refer to \cite{Guzhang2003} for details. Later, Gu improved the
restarted GMRES by augmenting the Krylov subspace with harmonic Ritz vectors for
Eq. (\ref{shsystem}) \cite{Gu2007}. By deflating eigenvalues for matrices that have
a few small eigenvalues, Darnell, Morgan and Wilcox \cite{Darnell08} presented an
improved GMRES method with deflated restarting to accelerate the convergence. Gu,
Zhou and Lin from another aspect of enhancing the convergence speed, proposed a
flexible preconditioned Arnoldi method that needs to exactly solve a linear system
with the coefficient matrix $A+\sigma_kI$ at the $k$-th iteration, where $\sigma_k$
is the precondition reference value that draws near $\alpha_j$. They also showed that
their proposed method is greatly faster than the traditional preconditioning strategies
\cite{GuZhoLin2007}. Saibaba, Bakhos and Kitanidis have further extended the flexible
preconditioning idea for solving generalized shifted linear systems arising from
oscillatory hydraulic tomography \cite{Saibaba2013}. Sun, Huang and Jing et al. 
\cite{SunHJC2018,Sun2018x} promoted the block version of GMRES method with deflated 
restarting for solving linear systems with multiple shifts and multiple right-hand 
sides. For other related methods, refer oneself to some studies in \cite{Frommer03,
Soodhalter14,GuHuaMen2014,JinHua2009,WuWanJin2012,Simon2010,Du2015,GijzenSlei15,Freund93,
Gu2016JCAM} and references therein.

As a cheaper implementation of GMRES, the Simpler GMRES algorithm (SGMRES) is
another famous Krylov subspace method \cite{WalZho1994}. It runs the Arnoldi
process begin with $Ar_0$ instead of $r_0$, where $r_0=b-Ax_0$. At each iteration,
it only requires to solve an upper-triangular least-squares problem rather than
an upper Hessenberg least-squares problem of GMRES, thus the SGMRES solver often
spends less computational cost. Recently, Jing, Yuan and Huang applied the SGMRES
and its stable variant: adaptive SGMRES (Ad-SGMRES) to solve the shifted system
(\ref{shsystem}) \cite{JinYuaHua2016}. For dealing with the non-colinearity of
$r_m$ and $r_m(\alpha_j)$, Jing, Yuan and Huang provided a remedy by forcing $
r_m(\alpha_j)\perp A\mathcal{K}_m(A,r_0)$. Besides to this advanced point, at
each iteration step, from the non-converged systems, they took the linear system
with the maximum residual norm as the seed system of the restart iteration.

However, in each cycle of the restarted methods, the convergence will slow down,
since the dimension of the Krylov subspace is limited \cite{GuZhoLin2007,LinBaoWu2012,
BooBhu2004,BooTanDoo2010,ZhoWuChe2015}. Especially for the problem with $A+\alpha_jI$
having small eigenvalues (in modulas). The main reason is that at each cycle, the
Krylov subspace does not contain good approximations of the eigenvectors corresponding
to such small eigenvalues. These make the thick-restarting and preconditioning techniques
beneficial for solving Eq. (\ref{shsystem}). Unfortunately, as far as we know, unlike
the shifted GMRES, there are not so many improved strategies applied to accelerate SGMRES
for solving shifted linear systems (\ref{shsystem}). Thus, in this paper, we will first
apply the flexible preconditioning technique \cite{SzyVog2001} to the Ad-SGMRES for solving
shifted linear systems (\ref{shsystem}), then consider restarting the new algorithm with
the deflated restarting strategy introduced in \cite{BooBhu2004,BooTanDoo2010}. The flexible preconditioning technique we used in this paper is the inexact preconditioning \cite{SimSzyld03} instead of exact which used in \cite{GuZhoLin2007}. The details will be located in Section \ref{sec2}.

The rest of this paper is organized as follows. In Section \ref{sec2}, we first give a brief description
of the adaptive Simpler GMRES method (Ad-SGMRES), then present two variants of Ad-SGMRES for shifted linear
system (\ref{shsystem}). Numerical examples in Section \ref{sec3} will illustrate the effectiveness of
the proposed algorithms. In Section \ref{sec4}, the paper closes with some conclusions.

\section{A flexible and adaptive Simpler GMRES algorithm with deflated restarting for shifted linear systems}
\label{sec2}
\setcounter{equation}{0}
\quad\
In this section, applying the flexible preconditioning technique \cite{SzyVog2001,Saibaba2013},
we first derive a flexible adaptive Simpler GMRES algorithm (FAd-SGMRES-Sh) for solving shifted
linear systems (\ref{shsystem}) simultaneously. Then based on it, we thick-restart the new algorithm
by using the deflated restarting strategy \cite{BooBhu2004,BooTanDoo2010,Meng2014}. Hence, a
flexible and adaptive Simpler GMRES algorithm with deflated restarting (FAd-SGMRES-DR-Sh) will
be achieved for solving Eq. (\ref{shsystem}).

Before giving the new algorithms, we will first briefly review the adaptive Simpler GMRES method.
By introducing a threshold parameter $\nu\in[0,1]$, Jir\'{a}nek and Rozlo\u{z}n\'{\i}k proposed
the adaptive Simpler GMRES (Ad-SGMRES) \cite{JirRoz2010}, which is more stable than the Simpler
GMRES, for solving the linear system $Ax = b$. The following algorithm is just the practical implementation of Ad-SGMRES.

\begin{algorithm}
{\bf The adaptive Simpler GMRES (Ad-SGMRES)}\\
{\bf 1.}~Given the initial guess $x_0$, a tolerance tol, a threshold parameter $\nu\in[0,1]$, let $m$ the maximal dimension of the solving subspace, $r_0=b-Ax_0$;\\
{\bf 2.}~For $k=1,\cdots,m$, do
\begin{flushleft}
~~~~~~~~~~~~(1)~~$z_k=\left\{\begin{array}{cc}
                                          r_0/\|r_0\|_2, & if~~ k=1, \\
                                          r_{k-1}/\|r_{k-1}\|_2, &~~~~~~~~~~~~~~~~~~ if~~ k>1,~~ and ~~\|r_{k-1}\|_2\leq \nu\|r_{k-2}\|_2, \\
                                          v_{k-1}, & otherwise.
                                        \end{array}\right.$\\
~~~~~~~~~~~~(2)~~$v_k=Az_k$,\\
~~~~~~~~~~~~(3)~~for $i=1,\cdots,k-1$\\
~~~~~~~~~~~~~~~~~~~~~$u_{ik}=v_i^Hv_k$,~~$v_k=v_k-u_{ik}v_i$.\\
~~~~~~~~~~~~~~~~~~end\\
~~~~~~~~~~~~(4)~~$u_{kk}=\|v_k\|_2$, $v_k=v_k/\|v_k\|_2$.\\
~~~~~~~~~~~~(5)~~$\xi_k=v_k^Hr_{k-1}$, $r_k=r_{k-1}-v_k\xi_k$, if $\|r_k\|_2\leq tol$, then go to {\bf Setp 3}.\\
~~~~~~~end
\end{flushleft}
{\bf 3.}~Let $k$ be the final iteration number of {\bf Step 2}, solve: $y_k=U_k^{-1}[\xi_{1},\cdots,\xi_k]^H$. Set $x_k=x_0+Z_ky_k$.
\label{Ad-SGMRES}
\end{algorithm}
In Algorithm \ref{Ad-SGMRES}, the definitions of $U_k$ and $V_k$ can be found in the next section.

\subsection{Flexible preconditioning}

\setcounter{equation}{0}
\quad\
Suppose $r_0=b-Ax_0\neq0$, where $x_0$ is the initial guess. At $k$-th iteration of Ad-SGMRES
(stated in Algorithm \ref{Ad-SGMRES}) for solving the seed system $Ax=b$, we have
\begin{equation}\label{SGMRESrel}
AZ_k=V_kU_k,
\end{equation}
where $Z_k=[z_1,\cdots,z_k]\in\mathbb{C}^{n\times k}$ is the basis of $\mathcal{K}_k(A,r_0)$, $V_k=[v_1,\cdots,v_k]\in\mathbb{C}^{n\times k}$ is the orthogonal basis of $A\mathcal{K}_k(A,
r_0)$, $U_k=[u_{ij}]\in\mathbb{C}^{k\times k},~i,j=1,\cdots,k$ is upper triangular, so $U_k$
is non-singular because the coefficient matrix $A$ is non-singular.

In \cite{GuZhoLin2007}, Gu, Zhou and Lin proposed a flexible preconditioning strategy for
GMRES that it is needed to exactly solve a linear system with the coefficient matrix $A +
\sigma_kI$ at the $k$-th iteration, and it will cost a lot of time especially for large size
problems. In this section, we will use the inexact flexible preconditioning \cite{SimSzyld03,
SzyVog2001,Sadd1993} instead of exact. It is known that the traditional right preconditioning
is applied to solve a modified system such as $AM^{-1}(Mx)=b$, where $AM^{-1}$ is well conditioned.
The inexact flexible preconditioning is actually a modification to the right preconditioning,
i.e., $M_k$ replaces $M$, so that inexact solver can be used. Based on such ideas, at each $k$-th
iteration, we set $w_k=M_k^{-1}z_k$, where $M_k$ is a variable preconditioner. Denote $W_k=[w_1,
\cdots,w_k]$, obviously, the columns of $W_k$ may not span a Krylov subspace. For the absence of
misunderstanding, we still use notions $V_k$ and $U_k$. The relation (\ref{SGMRESrel}) can be
rewritten in the following matrix equation:
\begin{equation}\label{FSGMRESrel}
AW_k=V_kU_k.
\end{equation}
For seed system, we seek the approximate solution $x_k=x_0+W_ky_k$ in the affine subspace $x_0
+ span\{W_k\}$, $y_k\in\mathbb{C}^k$ is a vector to be determined. Meanwhile, we seek the approximate solution $x_k(\alpha_j)=x_0(\alpha_j)+W_ky_k(\alpha_j)$ in the affine subspace $x_0(\alpha_j)+span\{W_k\}$ for add systems, where $y_k(\alpha_j)\in\mathbb{C}^k$ is a vector to be determined. For the add systems, we have
\begin{equation*}
\begin{split}
(A+\alpha_jI)W_k & = AW_k+\alpha_jW_k \\
& = V_kU_k+\alpha_jW_k.
\end{split}
\end{equation*}
Since for $W_k$ cannot be expressed by $V_k$, therefore, similar as in SGMRES \cite{WalZho1994}, there exists no $U_k(\alpha_j)$ for the add systems to keep a similar relation to (\ref{FSGMRESrel}). Hence, it is impossible to force the residual vectors $r_k(\alpha_j)$ to be colinear to $r_k$.

For the seed system $Ax=b$, since the orthogonal condition is $r_k\perp span\{AW_k\}$, i.e., $r_k\perp span\{V_k\}$, then using (\ref{FSGMRESrel}), we get
\begin{equation}
\begin{split}
0 & = V_k^H(b-Ax_k)\\
& = V_k^H(r_0-AW_ky_k)\\
& =V_k^Hr_0-U_ky_k,
\label{seedyk1}
\end{split}
\end{equation}
and
\begin{equation}
\begin{split}
r_k & = b-Ax_k \\
& = r_0-V_kU_ky_k \\
& = r_0-V_kV_k^Hr_0 \\
& =r_{k-1}-v_k\xi_k,
\end{split}
\label{seedres}
\end{equation}
where $\xi_k=v_k^Hr_0=v_k^Hr_{k-1}$. Thus (\ref{seedyk1}) can be rewritten as
\begin{equation}\label{seedyk2}
[\xi_1,\cdots,\xi_k]^H=U_ky_k.
\end{equation}
Similar to the strategy in \cite{JinYuaHua2016}, for the add systems, we require the residual vector $r_k(\alpha_j)=b-(A+\alpha_jI)x_k(\alpha_j)$ being orthogonal to $span\{AW_k\}$, together with (\ref{FSGMRESrel}), we have
\begin{equation}\label{addyk}
\begin{split}
0&=V_k^H[b-(A+\alpha_jI)x_k(\alpha_j)]\\
&=V_k^H(r_0(\alpha_j)-(AW_k+\alpha_jW_k)y_k(\alpha_j))\\
&=V_k^Hr_0(\alpha_j)-(U_k+\alpha_jV_k^HW_k)y_k(\alpha_j).
\end{split}
\end{equation}
Thus, after solving (\ref{seedyk2}) and (\ref{addyk}) to obtain $y_k$ and
$y_k(\alpha_k)$, the approximate solution of (\ref{shsystem}) is immediately
accessed, and then
\begin{equation}\label{addres}
r_k(\alpha_j)=r_0(\alpha_j)-(AW_k+\alpha_jW_k)y_k(\alpha_j)=r_0(\alpha_j)-(V_kU_k+\alpha_jW_k)y_k(\alpha_j).
\end{equation}

With the same seed system selection strategy in \cite{JinYuaHua2016,LiLiu2013}, we
summarize our flexible and adaptive Simpler GMRES for solving shifted linear systems
(FAd-SGMRES-Sh) in Algorithm 2. If the $\alpha_1$ in seed system is not zero, we can
reset
\begin{eqnarray*}
A&\doteq&A-\alpha_1I,\\
\alpha_j&\doteq&\alpha_j-\alpha_1,
\end{eqnarray*}
thus we take $\alpha_1=0$ as default.

\begin{algorithm}
{\bf A flexible and adaptive Simpler GMRES for shifted linear systems (FAd-SGMRES-Sh)}\\
{\bf 1.}~Start: Given the initial guess $x_0(\alpha_j)$, a tolerance tol, a threshold parameter $\nu\in[0,1]$, let $m$ the maximal dimension of the solving subspace, $r_0(\alpha_j)=b-Ax_0(\alpha_j)$;\\
{\bf 2.}~Select seed system: At the first iteration (after the second iteration), for all systems (for non-converged systems), find $ss\in\{1,\cdots,s\}$, where $s$ is adjusted by the number of non-converged systems, such that
\begin{equation*}
\|r_0(\alpha_{ss})\|_2=\max_{1\leq j\leq s}\|r_0(\alpha_j)\|_2.
\end{equation*}
Re-order $r_0(\alpha_1),\cdots,r_0(\alpha_s)$, so that the residual of the seed system is placed in the first place. Thus, after re-ordering, $ss=1$;\\
{\bf 3.}~Iterate: for $k=1,\cdots,m$, do
\begin{flushleft}
~~~~~~~~~~~~~~~~~~(1)~~$z_k=\left\{\begin{array}{cc}
                                          r_0/\|r_0\|_2, & if~~ k=1, \\
                                          r_{k-1}/\|r_{k-1}\|_2, &~~~~~~~~~~~~~~~~~~ if~~ k>1,~~ and ~~\|r_{k-1}\|_2\leq \nu\|r_{k-2}\|_2, \\
                                          v_{k-1}, & otherwise.
                                        \end{array}\right.$\\
~~~~~~~~~~~~~~~~~~(2)~~$w_k=M_k^{-1}z_k$,\\
~~~~~~~~~~~~~~~~~~(3)~~$v_k=Aw_k$,\\
~~~~~~~~~~~~~~~~~~(4)~~for $i=1,\cdots,k-1$\\
~~~~~~~~~~~~~~~~~~~~~~~~~~~$u_{ik}=v_i^Hv_k$,~~$v_k=v_k-u_{ik}v_i$.\\
~~~~~~~~~~~~~~~~~~~~~~~~end\\
~~~~~~~~~~~~~~~~~~(5)~~$u_{kk}=\|v_k\|_2$, $v_k=v_k/u_{kk}$.\\
~~~~~~~~~~~~~~~~~~(6)~~$\xi_k=v_k^Hr_{k-1}$, $r_k=r_{k-1}-v_k\xi_k$, if $\|r_k\|_2\leq tol$, then go to {\bf Setp 4}.\\
~~~~~~~~~~~~~end
\end{flushleft}
{\bf 4.}~Let $k$ be the final iteration number of {\bf Step 3}.\\
{\bf For seed system}, solve (\ref{seedyk2});\\
{\bf For add systems}, $j=2,\cdots,s$, solve (\ref{addyk}), and update $r_k(\alpha_j)$ using (\ref{addres});\\
{\bf 5.}~Set $x_k(\alpha_j)=x_0(\alpha_j)+W_ky_k(\alpha_j)$, $j=1,\cdots,s$. For the non-converged systems, reset $r_0(\alpha_j)=r_k(\alpha_j)$, $x_0(\alpha_j)=x_k(\alpha_j)$, $j=1,\cdots,s$, go to {\bf step 2}.
\label{FAd-SGMRES-Sh}
\end{algorithm}

Some remarks of the implementation details for FAd-SGMRES-Sh are as follows.
\begin{rem}
In {\bf Step 3}, $M_k$ is the flexible preconditioner in the $k$-th step. To
get the effect of preconditioning, $M_k$ is usually selected to be the matrix
near $A$. In our algorithm, we choose to solve $Aw_k=z_k$ inexactly for the
process $w_k=M_k^{-1}z_k$. There are many choices of inexact solvers, such as
ILU \cite{Saa2003}, IHSS \cite{BaiGolNg2003}, IGMRES \cite{Saa2003}, ISOR \cite{
Saa2003}, IQR \cite{BaiDufYin2009}, and so on. In numerical examples section,
we select IGMRES with 10 iterations as the preconditioner.
\end{rem}

\begin{rem}
In {\bf Step 4}, for add systems, the matrix $U_k+\alpha_jV_k^HW_k$ is generally
not upper triangular. Because we usually choose a small value $m\ll n$, such as
20, thus for the solving step $V_k^Hr_0(\alpha_j)=(U_k+\alpha_jV_k^HW_k)y_k(\alpha_j)$,
the MATLAB code ``$\setminus$" can be directly used to get $y_k(\alpha_j)$. In addition,
from (\ref{addres}), we can see the update of the residual vectors will also cost
some time. Consequently, for solving add systems, similar to SGMRES \cite{JinYuaHua2016},
FAd-SGMRES-Sh may not faster than GMRES \cite{Frommer03}. But fortunately, for seed
system, due to without solving an upper Hessenberg least-square problem, and with
inexact preconditioning, FAd-SGMRES-Sh is much faster than SGMRES, GMRES and FGMRES
\cite{GuZhoLin2007}, especially for large-scale problems. Numerical experiments will
illustrate the effect later.
\end{rem}
%%%%%%%%%%%%%%%%%%%%%%%%%%%%%%%%%%%%%%%%
\subsection{Thick-restarting}
\quad\
Actually, some inexact preconditioned systems may still encounter the issues with
small eigenvalues, thus it is necessary to consider to restart Algorithm \ref{FAd-SGMRES-Sh} with
the deflated restarting strategy \cite{GauGutLieNab2013,Mor2002,Meng2014}. Our aim
is to improve the convergence of FAd-SGMRES-Sh by using the spectral information of
the preconditioned seed system at restart. There are two keys involved. The first
is how to compute the spectral information at each restart. The second is how to
apply these information with a low computation cost at restart.

In fact, we use the harmonic Ritz value information of the seed system $Ax=b$ at 
each restart. That is required, after one cycle, the harmonic Ritz pair $(\lambda_i, 
q_i\equiv W_mg_i)$ of $A$ in $span\{W_m\}$ and orthogonal to $span\{AW_m\}$ satisfying \cite{PaiParVor1995}:
\begin{equation*}
AW_mg_i-\lambda_iW_mg_i\perp span\{AW_m\}\Leftrightarrow (V_m U_m)^H(AW_mg_i-\lambda_iW_mg_i)=0.
\end{equation*}
From (\ref{FSGMRESrel}), and $U_m$ non-singular, the above equation is equivalent to
\begin{equation}\label{FAdSGMRESDReig}
U_mg_i=\lambda_iV_m^HW_mg_i.
\end{equation}
Consequently, the harmonic Ritz pairs can be calculated at each iteration of FAd-SGMRES-Sh. Let $(\lambda_i,g_i), i=1,\cdots,e(e\leq m)$ are the eigenpairs of the reduced generalized eigenvalues problem (\ref{FAdSGMRESDReig}). Set $G_e=[g_1,\cdots,g_e]$, suppose that $P_eL_e=G_e$ is the QR decomposition of $G_e$, where matrix $P_e=[p_1,\cdots,p_e]\in\mathbb{C}^{k\times e}$ is orthogonal. Postmultiplying (\ref{FSGMRESrel}) by $P_e$ yields
\begin{equation}\label{FAdSGMRESdrrel1}
AW_mP_e=V_mU_mP_e.
\end{equation}
Let $U_mP_e=\widehat{P}_eU_e^{new}$ be the QR decomposition, then from (\ref{FAdSGMRESdrrel1}) we have
\begin{equation*}
AW_mP_e=V_m\widehat{P}_eU_e^{new}.
\end{equation*}
Define $W_e^{new}=W_mP_e$ and $V_e^{new}=V_m\widehat{P}_e$, then we obtain
\begin{equation*}
AW_e^{new}=V_e^{new}U_e^{new},
\end{equation*}
where $V_e^{new}\in\mathbb{C}^{n\times e}$ is orthogonal, $U_e^{new}\in\mathbb{C}^{e\times e}$ is upper triangular. Let $W_e=W_e^{new}$, $V_e=V_e^{new}$ and $U_e=U_e^{new}$. To establish the equation (\ref{FSGMRESrel}) for the current cycle, the flexible and adaptive Simpler GMRES with deflated restarting executes the remaining ($m-e$) steps with $w_i=M_i^{-1}z_i (e+1\leq i\leq m)$ %as the initial vector,
where $M_i$ is the flexible preconditioner and
\begin{equation*}
z_i=\left\{
  \begin{array}{cc}
     r_e/\|r_e\|_2, & if~ i=e+1, \\
     r_{i-1}/\|r_{i-1}\|_2, & ~~~~if~  i>e+1~ and~ \|r_{i-1}\|_2\leq\nu\|r_{i-2}\|_2, \\
     v_{i-1},& otherwise.
   \end{array}\right.
\end{equation*}

After each cycle of the new algorithm, we restart the algorithm by setting $x_0^{new}(\alpha_j)=x_m(\alpha_j)$ and $r_0^{new}(\alpha_j)=r_m(\alpha_j)$. We use the symbols such as $x_m^{new}(\alpha_j)$, $r_m^{new}(\alpha_j)$, $W_m^{new}$, $V_m^{new}$ and $U_m^{new}$ for current cycle to distinguish the ones from the last cycle.

For the seed system, after one cycle of FAd-SGMRES-Sh, from (\ref{seedres}), we have
\begin{equation*}
r_0^{new}=r_m=r_0-V_mV_m^Hr_0,
\end{equation*}
and
\begin{equation*}
r_e^{new}=r_0^{new}-V_e^{new}(V_e^{new})^Hr_0^{new}.
\end{equation*}
Note that
\begin{equation*}
(V_e^{new})^Hr_0^{new}=\widehat{P}_e^HV_m^H(r_0-V_mV_m^Hr_0)=0.
\end{equation*}
Thus
\begin{equation*}
r_e^{new}=r_0^{new},~~~~\xi_i^{new}=(v_i^{new})^Hr_0^{new}=0,~ ~~~i=1,\cdots,e,
\end{equation*}
then from (\ref{seedyk1}) and (\ref{seedres}), we need to solve \begin{equation}\label{FAdSGMRESDRseedy}
U_m^{new}y_m^{new}=[0,\cdots,0,\xi_{e+1}^{new},\cdots,\xi_m^{new}]^T,
\end{equation}
where $\xi_i^{new}=(v_i^{new})^Hr_0^{new}=(v_i^{new})^Hr_{i-1}^{new}$, $i=e+1,\cdots,m$, and update
\begin{equation}\label{FAdSGMRESseedr}
r_i^{new}=r_{i-1}^{new}-v_i^{new}\xi_i^{new}.
\end{equation}

For add systems, from (\ref{addyk}) we can get
\begin{equation*}
(V_e^{new})^Hr_0(\alpha_j)^{new}=\widehat{P}_e^HV_m^Hr_m(\alpha_j)=0,
\end{equation*}
thus,
\begin{equation*}
(V_m^{new})^Hr_0(\alpha_j)^{new}=[0,\cdots,0,\xi_{e+1}(\alpha_j)^{new},\cdots,\xi_m(\alpha_j)^{new}]^T,
\end{equation*}
where $\xi_i(\alpha_j)^{new}=(v_i^{new})^Hr_0(\alpha_j)^{new}$, $i=e+1,\cdots,m$. Consequently, from (\ref{addyk}), we need to solve
\begin{equation}\label{FDraddyk}
[0,\cdots,0,\xi_{e+1}(\alpha_j)^{new},\cdots,\xi_m(\alpha_j)^{new}]^T=(U_k^{new}+\alpha_j(V_k^{new})^HW_k^{new})y_k(\alpha_j)^{new},
\end{equation}
and we still exploit (\ref{addres}) to update the residual vector. Now it
is ready to present the main algorithm of this paper.

\begin{algorithm}
{\bf A flexible and adaptive Simpler GMRES with deflated restarting for shifted linear systems (FAd-SGMRES-DR-Sh)}\\
{\bf 1.}~Start: Given the initial guess $x_0(\alpha_j)$, an integer $e$, a tolerance tol, a threshold parameter $\nu\in[0,1]$, let $m$ the maximal dimension of the solving subspace,  $r_0(\alpha_j)=b-x_0(\alpha_j)$;\\
{\bf 2.}~Select seed system: At the first iteration (after the second iteration), for all systems (for non-converged systems), find $ss\in\{1,\cdots,s\}$, where $s$ is adjusted by the number of non-converged systems, such that
\begin{equation*}
\|r_0(\alpha_{ss})\|_2=\max_{1\leq j\leq s}\|r_0(\alpha_j)\|_2.
\end{equation*}
Re-order $r_0(\alpha_1),\cdots,r_0(\alpha_s)$, so that the residual of the seed system is placed in the first place. Thus, after re-ordering, $ss=1$;\\
{\bf 3.}~Apply one cycle of FAd-SGMRES-Sh to the seed system $Ax=b$, generate $W_m$, $V_m$, $U_m$, $x_m$, and $r_m$;\\
{\bf 4.}~Compute the eigenvalues and eigenvectors of the generalized eigenvalue problem (\ref{FAdSGMRESDReig}) by using the QZ algorithm. Let $g_1,\cdots,g_e$ be the eigenvectors corresponding to the $e$ smallest eigenvalues of (\ref{FAdSGMRESDReig}). Set $G_e=[g_1,\cdots,g_e]$, and compute the QR decompositions of $G_e$ and $U_mP_e$: $G_e=P_eL_e$, $U_mP_e=\widehat{P}_eU_e^{new}$. Let $W_e^{new}=W_mP_e$ and $V_e^{new}=V_m\widehat{P}_e$.\\
{\bf 5.}~Let $W_e=W_e^{new}$, $V_e=V_e^{new}$, $U_e=U_e^{new}$, and $x_0=x_m$, $r_0=r_m$, $r_e=r_0$;\\
{\bf 6.}~Iterate: for $k=e+1,\cdots,m$, do
\begin{flushleft}
~~~~~~~~~~~~~~~~~~(1)~~$z_k=\left\{\begin{array}{cc}
                                          r_e/\|r_e\|_2, & if~~ k=e+1, \\
                                          r_{k-1}/\|r_{k-1}\|_2, &~~~~~~~~~~~~~~~~~~ if~~ k>e+1,~~ and ~~\|r_{k-1}\|_2\leq \nu\|r_{k-2}\|_2, \\
                                          v_{k-1}, & otherwise.
                                        \end{array}\right.$\\
~~~~~~~~~~~~~~~~~~(2)~~$w_k=M_k^{-1}z_k$,\\
~~~~~~~~~~~~~~~~~~(3)~~$v_k=Aw_k$,\\
~~~~~~~~~~~~~~~~~~(4)~~for $i=1,\cdots,k-1$\\
~~~~~~~~~~~~~~~~~~~~~~~~~~~$u_{ik}=v_i^Hv_k$,~~$v_k=v_k-u_{ik}v_k$.\\
~~~~~~~~~~~~~~~~~~~~~~~~end\\
~~~~~~~~~~~~~~~~~~(5)~~$u_{kk}=\|v_k\|_2$, $v_k=v_k/\|v_k\|_2$.\\
~~~~~~~~~~~~~~~~~~(6)~~$\xi_k=v_k^Hr_{k-1}$, $r_k=r_{k-1}-v_k\xi_k$, if $\|r_k\|_2\leq tol$, then go to {\bf Setp 7}.\\
~~~~~~~~~~~~~end
\end{flushleft}
{\bf 7.}~Let $k$ be the final iteration number of {\bf Step 6}.\\
{\bf For seed system}, solve (\ref{FAdSGMRESDRseedy});\\
{\bf For add systems}, $j=2,\cdots,s$, solve (\ref{addyk}), and
update $r_k(\alpha_j)$ using (\ref{addres});\\
{\bf 8.}~Set $x_k(\alpha_j)=x_0(\alpha_j)+W_ky_k(\alpha_j)$, $j=1,\cdots,s$. For the non-converged systems, reset $r_0(\alpha_j)=r_k(\alpha_j)$, $x_0(\alpha_j)=x_k(\alpha_j)$, $j=1,\cdots,s$, go to {\bf step 2}.
\label{FAd-SGMRES-DR-Sh}
\end{algorithm}

In the end of this section, it is meaningful to evaluate the computational costs in a generic cycle of GMRES-Sh,
Ad-SGMRES-Sh, FAd-SGMRES-Sh and FAd-SGMRES-DR-Sh, where the detail pseudo-codes of GMRES-Sh and Ad-SGMRES-Sh are be found in \cite{JinYuaHua2016}.
The comparisons are presented in Table \ref{tab1} and Table \ref{tab2}. Here, we denote ``mv" the number of matrix-vector products.
``$op_{M_k}$" denots the number of the preconditioning process $M_k^{-1}z_k$ in one cycle, ``vector updates" denotes
the number of vectors that need to be updated in one cycle. We also write down the number of generalized eigenvalue problems
by ``G-p" in one cycle.
%\medskip
%{
\begin{table}[!htbp]
\centering
\caption{Main computational costs per cycle for GMRES-Sh, Ad-SGMRES-Sh and FAd-SGMRES-Sh}
\vspace{0.5mm}
\begin{tabular}{cccc}
\hline
&GMRES-Sh&Ad-SGMRES-Sh&FAd-SGMRES-Sh\\
\hline
mv&$m$&$m$&$m$\\
dot products&$m(\sum\limits_{k=1}^{m}(k-1)+1)$&$m(\sum\limits_{k=1}^{m}(k-1)+1+s)$&$m(\sum\limits_{k=1}^{m}(k-1)+1+s)$\\
saxpy&$m(\sum\limits_{k=1}^{m}(k-1)+1)+m+s$&$m(\sum\limits_{k=1}^{m}(k-1)+1)+2s$&$m(\sum\limits_{k=1}^{m}(k-1)+1)+2s$\\
$op_{M_k}$&0&0&$m$\\
%Matr-seed&$\overline{H}_m$&$U_m$&$U_m$&$U_m$&$U_m$\\
vector updates&$m+s+1$&$2m+2s$&$2m+2s$\\
%Matr-add&$[\overline{H}_m+\alpha_j,z]$&$U_m+\alpha_jV_m^HZ_m$&$U_m+\alpha_jV_m^HZ_m$&$U_m+\alpha_jV_m^HZ_m$&$U_m+\alpha_jV_m^HZ_m$\\
G-p&0&0&0\\
\hline
\end{tabular}
\label{tab1}
\end{table}

\begin{table}[!htpb]
\centering
\caption{Main computational costs per cycle for the 1st cycle and the other cycle of FAd-SGMRES-DR-Sh}
\vspace{0.5mm}
\begin{tabular}{ccc}
\hline
&FAd-SGMRES-DR-Sh&FAd-SGMRES-DR-Sh\\
&(1st cycle)&(other cycle)\\
\hline
mv&$m$&$m-e$\\
dot products&$m(\sum\limits_{k=1}^{m}(k-1)+1+s)$&$(m-e)(\sum\limits_{k=1}^{m}(k-1)+1+s)$\\
saxpy&$m(\sum\limits_{k=1}^{m}(k-1)+1)+2s$&$(m-e)(\sum\limits_{k=1}^{m}(k-1)+1)+2s$\\
$op_{M_k}$&$m$&$m-e$\\
%Matr-seed&$\overline{H}_m$&$U_m$&$U_m$&$U_m$&$U_m$\\
vector updates&$2m+2s$&$2m+2s$\\
%Matr-add&$[\overline{H}_m+\alpha_j,z]$&$U_m+\alpha_jV_m^HZ_m$&$U_m+\alpha_jV_m^HZ_m$&$U_m+\alpha_jV_m^HZ_m$&$U_m+\alpha_jV_m^HZ_m$\\
G-p&1&1\\
\hline
\end{tabular}
\label{tab2}
\end{table}

\section{Numerical results}
\label{sec3}
\quad\
In this section, numerical comparisons are made for GMRES-Sh \cite{Frommer03},
Ad-SGMRES-Sh \cite{JinYuaHua2016}, FGMRES-Sh \cite{GuZhoLin2007}, GMRES-DR-Sh
\cite{Mor2002}, FAd-SGMRES-Sh and FAd-SGMRES-Dr-Sh according to the number of
outer matrix-vector products (referred to as $mv$), and the elapsed CPU time
in seconds (referred to as $cpu$). %Here, the number of inner matrix-vector products means the products cost in the preconditioning process, i.e., the {\bf Step 3 (2)} in Algorithm 2 and the {\bf Step 6 (2)} in Algorithm 3. Because there is no preconditioning in Ad-SGMRES-Sh, then $imv=0$.
We set the
stopping criterion as
\begin{equation*}
\frac{\|b-(A+\alpha_jI)x_k(\alpha_j)\|_2}{\|b\|_2} < {\tt 1e-6},\quad~ j=1,2,\cdots,s.
\end{equation*}
The bold values
in the following tables indicate the fastest in the terms of $cpu$. The numerical
results are obtained by using MATLAB R2014a (64bit) on an PC-Intel Core i5-6200U,
CPU 2.4 GHz, 8 GB RAM with machine epsilon $10^{-16}$ in double precision floating
point arithmetic.

{\bf Example 3.1} We consider the same matrices used in \cite{JinYuaHua2016}.
These matrices are from the University of Florida Sparse Matrix Collection and the
Example 1 in \cite{Mor2005}. Table \ref{tab3} lists the matrices with their information.
Here bidiag1 and bidiag2 are bidiagonal matrices with super-diagonal entries being all one.
The diagonal elements of bidiag1 are $0.1, 1,2,3,\cdots,999$, and the ones of bidiag2
are $1,2,3,\cdots,1000$. All the initial vectors are zero in all examples. The right-hand
side $b$ is generated by the MATLAB code $randn(n,1)$, where $n$ is the dimension of $A$.
The shift parameters are $\alpha=0,0.4,2$. For FAd-SGMRES-Sh and FAd-SGMRES-DR-Sh, the
flexible preconditioner is chosen as running 10 steps of the un-restarted GMRES algorithm
\cite{SaaSch1986}. The same strategy is used in Example 3.2. For FGMRES-Sh, we use LU decomposition to exactly solve $(A+\sigma_i)w=v$ in the preconditioning process. Similar as in \cite{GuZhoLin2007},  we select the same $\sigma_1=0.5$ in the first $m/2$ steps, in the last $m/2$ steps for the same $\sigma_2=1$. Thus, the LU decomposition of $A+\sigma_iI$ need to save for using in the first and last $m/2$ steps of each cycle. The same strategy is also used in Example 3.3.

In Table \ref{tab4}, we reported the $mv(cpu)$ of each algorithm for listed matrices with
size smaller than 1000, and the dimension of the approximate subspace in each cycle is set 
as $m=10$, $\mu=0.9$. For FAd-SGMRES-Dr-Sh, $e$ is the number of harmonic eigenvectors 
retained from the previous cycle. We compare two cases, i.e., $e=3,6$. In Table \ref{tab5}, 
for comparison, we set $m=20$ and $e=5,10,15$, with $\mu=0.9$, and the matrices size are 
all larger than 1000. In all tables, ``$\dag$" stands for the algorithm fails to converge 
even after using 10000 outer matrix-vector products.
\begin{table}[!htpb]
\centering
\caption{The test matrices used in Example 4.1}
\vspace{0.5mm}
\begin{tabular}{ccccc}
\hline
Matrix ID & Matrix name & Size   & Nonzeros & Problem domain               \\
\hline
1       & add20       & 2,395  & 13,151   & Circuit simulation           \\
2       & bidiag1     & 1,000  & 1,999    & Academic                     \\
3       & bidiag2     & 1,000  & 1,999    & Academic                     \\
4       & cdde1       & 961    & 4,681    & Computational fluid dynamics \\
5       & epb1        & 14,734 & 95,053   & Thermal                      \\
6       & sherman4    & 1,104  & 3,786    & Computational fluid dynamics \\
7       & wang1       & 2,903  & 19,093   & Semiconductor device         \\
8       & wang4       & 26,068 & 177,196  & Semiconductor device         \\
9       & young1c     & 841    & 4,089    & Acoustics                    \\
10      & young2c     & 841    & 4,089    & Acoustics                    \\
\hline
\end{tabular}
\label{tab3}
\end{table}

\begin{table}[t]
\centering
\caption{Convergence behaviors of the GMRES-Sh, Ad-SGMRES-Sh, FGMRES-Sh,
FAd-SGMRES-Sh, GMRES-DR-Sh and FAd-SGMRES-DR-Sh with ${\tt tol = 1e-6}$, $m=10$ and $\mu=0.9$}
\vspace{0.5mm}
\begin{tabular}{cccccc}
\hline
Method & \multicolumn{5}{c}{$mv(cpu),~m=10,~\mu=0.9$}\\
[-2pt]\cmidrule(l{0.7em}r{0.7em}){2-5}\\[-11pt]
  &bidiag1&bidiag2&cdde1&young1c&young2c\\
  \hline
  GMRES-Sh&4678(0.36)&513(0.06)&$\dag$&$\dag$&$\dag$\\
  Ad-SGMRES-Sh&4678(0.34)&513(0.06)&9569(1.00)&$\dag$&$\dag$\\
  FGMRES-Sh&7({\bf 0.02})&7({\bf 0.01})&118(0.07)&12({\bf 0.13})&11({\bf 0.12})\\
  FAd-SGMRES-Sh&54(0.04)&35(0.03)&21(0.06)&627(0.73)&615(0.72)\\
  GMRES-DR-Sh&\multirow{2}{*}{351(0.81)}&\multirow{2}{*}{258(0.05)}&\multirow{2}{*}{174(0.10)}&\multirow{2}{*}{$\dag$}&\multirow{2}{*}{$\dag$}\\
  $e=3$&&&&&\\
  GMRES-DR-Sh&\multirow{2}{*}{373(0.10)}&\multirow{2}{*}{240(0.06)}&\multirow{2}{*}{169(0.06)}&\multirow{2}{*}{$\dag$}&\multirow{2}{*}{$\dag$}\\
  $e=6$&&&&&\\
  FAd-SGMRES-DR-Sh&\multirow{2}{*}{39(0.02)}&\multirow{2}{*}{32(0.02)}&\multirow{2}{*}{19(0.06)}&\multirow{2}{*}{231( 0.34)}&\multirow{2}{*}{230(0.34)}\\
  $e=3$&&&&&\\
  FAd-SGMRES-DR-Sh&\multirow{2}{*}{41(0.03)}&\multirow{2}{*}{32(0.02)}&\multirow{2}{*}{19({\bf 0.04})}&\multirow{2}{*}{193( 0.25)}&\multirow{2}{*}{178(0.24)}\\
  $e=6$&&&&&\\
  \hline
\end{tabular}
\label{tab4}
\end{table}

\begin{table}
\centering
\caption{Convergence behaviors of the GMRES-Sh, Ad-SGMRES-Sh, FGMRES-Sh, 
FAd-SGMRES-Sh, GMRES-DR-Sh and FAd-SGMRES-DR-Sh with ${\tt tol = 1e-6}$,
$m=20$ and $\mu=0.9$}
\vspace{0.5mm}
\begin{tabular}{cccccc}
  \hline
  \multicolumn{6}{c}{$mv(cpu),~m=20,~\mu=0.9$}\\
  \cline{2-6}
  Method&add20&epb1&sherman4&wang1&wang2\\
  \hline
  GMRES-Sh&1231(0.32)&1300(1.17)&548(0.17)&1049(0.33)&$\dag$\\
  Ad-SGMRES-Sh&1231(0.27)&1310(1.70)&548(0.10)&894(0.24)&$\dag$\\
  FGMRES-Sh&635(11.39)&1099(9.74)&14(0.07)&295(1.38)&3161(268.17)\\
  FAd-SGMRES-Sh&55(0.10)&72(0.59)&23(0.05)&51(0.12)&148(2.30)\\
  GMRES-DR-Sh&\multirow{2}{*}{629(0.27)}&\multirow{2}{*}{601(1.89)}&\multirow{2}{*}{134(0.09)}&\multirow{2}{*}{473(0.25)}&\multirow{2}{*}{1162(7.05)}\\
  $e=5$&&&&&\\
  GMRES-DR-Sh&\multirow{2}{*}{$\dag$}&\multirow{2}{*}{591(2.30)}&\multirow{2}{*}{130(0.04)}&\multirow{2}{*}{496(0.23)}&\multirow{2}{*}{$\dag$}\\
  $e=10$&&&&&\\
  GMRES-DR-Sh&\multirow{2}{*}{$\dag$}&\multirow{2}{*}{$\dag$}&\multirow{2}{*}{136(0.06)}&\multirow{2}{*}{$\dag$}&\multirow{2}{*}{$\dag$}\\
  $e=15$&&&&&\\
  FAd-SGMRES-DR-Sh&\multirow{2}{*}{56(0.11)}&\multirow{2}{*}{63({\bf 0.58})}&\multirow{2}{*}{23(0.06)}&\multirow{2}{*}{44(0.11)}&\multirow{2}{*}{80({\bf 1.25)}}\\
  $e=5$&&&&&\\
  FAd-SGMRES-DR-Sh&\multirow{2}{*}{55({\bf 0.08})}&\multirow{2}{*}{63(0.59)}&\multirow{2}{*}{23(0.02)}&\multirow{2}{*}{{\bf 44(0.07)}}&\multirow{2}{*}{76(1.26)}\\
  $e=10$&&&&&\\
  FAd-SGMRES-DR-Sh&\multirow{2}{*}{56(0.09)}&\multirow{2}{*}{63(0.70)}&\multirow{2}{*}{{\bf 23(0.01)}}&\multirow{2}{*}{44(0.08)}&\multirow{2}{*}{79(1.50)}\\
  $e=15$&&&&&\\
  \hline
\end{tabular}
\label{tab5}
\end{table}
 %\begin{tabular*}{}{@{\extracolsep\fill}ccccc}
  %\hline
  %\multicolumn{5}{c}{$mv(cpu),~m=10,~\mu=0.9$}\\
  %\cline{2-5}
  %\multirow{2}{*}{Matrix}&\multirow{2}Ad-SGMRES-Sh}&\multirow{2}{*}{FAd-SGMRES-Sh}&\multicolumn{2}{c}{FAd-SGMRES-Dr-Sh}\\
  %\cline{4-5}
  %&&&$e=3$&$e=6$\\
  %\hline
  %add20&2799(0.406)&64(0.078)&55({\bf 0.068})&58(0.078)\\
  %bidiag1&$\dag$&56(0.047)&40({\bf 0.031})&43(0.047)\\
  %bidiag2&564(0.062)&36(0.032)&33({\bf 0.028})&34(0.031)\\
  %cdde1&945(0.094)&19(0.017)&18({\bf 0.015})&18(0.016)\\
  %epb1&1757(1.735)&72(0.547)&67({\bf 0.516})&67(0.531)\\
  %sherman4&2557(2.5)&84(0.641)&66({\bf 0.515})&67(0.594)\\
  %wang1&1770(0.296)&74(0.109)&47({\bf 0.078})&48(0.079)\\
  %wang4&2690(4.344)&164(2.203)&86({\bf 1.234})&83(1.297)\\
  %young1c&$\dag$&586(0.625)&233(0.266)&178({\bf 0.234})\\
  %young2c&$\dag$&621(0.64)&245(0.281)&172({\bf 0.218})\\
  %\hline
%\end{tabular*}}

%\medskip
As seen from Table \ref{tab4} and Table \ref{tab5}, for smaller matrices except for 
cdde1, FGMRES-Sh is the best solver among these algorithms, which is inseparable from
the exact solution of $(A+\sigma_k)w=v$ during the preconditioning process. But for 
the larger matrices, especially for wang4 whose size is 26068, the exact solving process 
of FGMRES-Sh obviously became a time-consuming obstacle, while FAd-SGMRES-DR-Sh performs best. It 
also can see for FAd-SGMRES-DR-Sh with different values $e$, in some examples, e.g., 
epb1 in Table \ref{tab5}, even the number $mv$ is smaller, but the elapsed CPU time 
is larger, this is because when using the harmonic Ritz value information, there needs 
to compute a generalized eigenvalue problem (\ref{FAdSGMRESDReig}) and sort these 
eigenvalues, thus if the eigenvectors number $e$ is larger, the elapsed CPU time for 
the previous procedure may be larger too. Thus, it is important to choose appropriate 
$m$ and $e$. For some matrices, such as bidiag2, cdde1, add20 and sherman4, we can 
see the number $mv$ of FAd-SGMRES-DR-Sh is not much less than FAd-SGMRES-Sh, even 
equal to each other, this is because after preconditioning, the small eigenvalues 
problems of these matrices are well controlled, thus the effect of deflated restarting 
is not obvious, whereas the other matrices are still need the deflated restarting. 
Consequently, for large and difficult problems, FAd-SGMRES-DR-Sh still performs better 
than the other mentioned algorithms.

\begin{figure}[!htpb]
\centering
\includegraphics[width=0.5\textwidth]{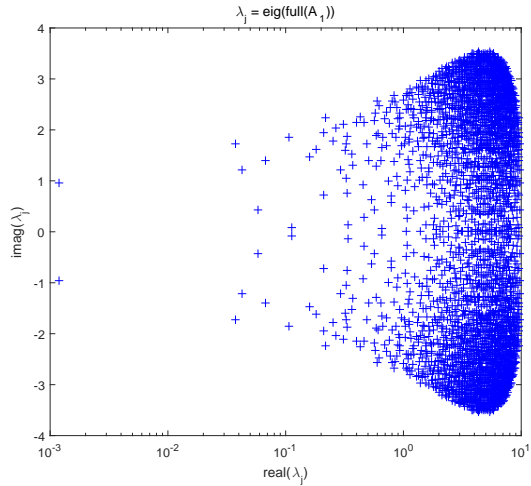}\\
\caption{The eigenvalues distribution of $A_1$.}
\label{fig21}
\end{figure}

{\bf Example 3.2} In this example, we apply our algorithms to solve quantum 
chromodynamics (QCD) problems with multiple shifts, which is one of the most 
time-consuming supercomputer applications. $D_i, 1\leq i\leq 14$ are denoted 
the complex matrices downloaded from Matrix Market\footnote{Refer to the website: 
\url{http://math.nist.gov/MatrixMarket/}.}. These $D_i$ are discretizations by 
the Dirac operator used in numerical simulation of quark behavior at different 
physical temperatures \cite{DarMorWil2004,Soodhalter14}. For each $D_i$, we take 
$A_i=(\frac{1}{k_c}+10^{-3})I-D_i$ as the base matrix, where $k_c$ is the critical 
value such that for $\frac{1}{k_c}<\frac{1}{k}<\infty$, the matrix $\frac{1}{k}I-
D_i$ is {\it real-positive}. Table \ref{tab6} lists the matrices $D_i$ with their information. Moreover, the right-hand side $b=\texttt{ones(length(A),1)}$, and the initial guess in each example is zero vector. We take $[0.0001,0.0002,\ldots,0.0004, 0.001,0.002,\ldots,0.004,0.01,0.02,\ldots,0.04]$ as the set of shifted values $\alpha_j$. It is shown from Fig. \ref{fig21} that the eigenvalues of base matrix $A_1$ are in the right-half of the complex plane, but partially surround the origin \cite{Darnell08}.

For seed matrices $A_1-A_7$, we set $m=10$, $\mu=0.9$, and $e=3,6$. Table \ref{tab7} 
gives the results of the considered algorithms. Form Table \ref{tab7}, it can see 
that GMRES-DR-Sh does not converge for each matrix, and FGMRES-Sh costs too much time, 
whereas FAd-SGMRES-DR-Sh performs best, this implies that after adding inexact 
preconditioning and then deflating the small eigenvalues can accelerate the 
convergence. In Table \ref{tab8}, we compares the other algorithms besides 
FGMRES-Sh and GMRES-DR-Sh, and we set $m=20$, $\mu=0.9$,
$e=5,10,15$ for seed matrices $A_8-A_{14}$. As seen from Table \ref{tab7} 
and Table \ref{tab8}, FAd-SGMRES-DR-Sh performs better than the other algorithms 
for most examples with deflating the small eigenvalues (in modulas). It is 
also known that the appropriate choice of $m$ and $e$ is important for 
FAd-SGMRES-DR-Sh, which will be subject to further investigations in the future.

%{\small
\begin{table}[t]
\centering
\caption{The matrices $D_i$ used in Example 4.2}
\vspace{0.5mm}
\begin{tabular}{cccccc}
\hline
Matrix ID &Denotation& Matrix name &Size & Nonzeros & $k_c$\\
\hline
1  &$D_1$    &CONF5.0-00L4X4-1000 &3,072  &119,808   &0.20611\\
2  &$D_2$    &CONF5.0-00L4X4-1400 &3,072  &119,808   &0.20328\\
3  &$D_3$    &CONF5.0-00L4X4-1800 &3,072  &119,808   &0.20265\\
4  &$D_4$    &CONF5.0-00L4X4-2200 &3,072  &119,808   &0.20235\\
5  &$D_5$    &CONF5.0-00L4X4-2600 &3,072  &119,808   &0.21070\\
6  &$D_6$    &CONF6.0-00L4X4-2000 &3,072  &119,808   &0.17968\\
7  &$D_7$    &CONF6.0-00L4X4-3000 &3,072  &119,808   &0.16453\\
8  &$D_8$    &CONF5.4-00L8X8-0500 &49,152 &1,916,928 &0.17865\\
9  &$D_9$    &CONF5.4-00L8X8-1000 &49,152 &1,916,928 &0.17843\\
10 &$D_{10}$ &CONF5.4-00L8X8-1500 &49,152 &1,916,928 &0.17689\\
11 &$D_{11}$ &CONF5.4-00L8X8-2000 &49,152 &1,916,928 &0.17835\\
12 &$D_{12}$ &CONF6.0-00L8X8-2000 &49,152 &1,916,928 &0.15717\\
13 &$D_{13}$ &CONF6.0-00L8X8-3000 &49,152 &1,916,928 &0.15649\\
14 &$D_{14}$ &CONF6.0-00L8X8-8000 &49,152 &1,916,928 &0.15623\\
  \hline
\end{tabular}
\label{tab6}
\end{table}

\begin{table}[!htpb]\small
\centering
\caption{Convergence behaviors of the Ad-SGMRES-Sh, FAd-SGMRES-Sh and 
FAd-SGMRES-DR-Sh with $n=3072$, ${\tt tol = 1e-6}$, $m=10$ and $\mu=0.9$}
\vspace{0.5mm}
\begin{tabular}{cccccccc}
  \hline
  \multicolumn{8}{c}{$mv(cpu),~m=10,~\mu=0.9$}\\
  \cline{2-8}
  Method&A1&A2&A3&A4&A5&A6&A7\\
  \hline
  GMRES-Sh&812(0.56)&315(0.36)&634(0.57)&384(0.40)&564(0.53)&2357(1.78)&176(0.26)\\
  Ad-SGMRES-Sh&812(0.86)&315(0.36)&634(0.70)&384(0.44)&564(0.62)&2357(2.53)&176(0.22)\\
  FGMRES-Sh&6(21.67)&6(32.21)&5(21.91)&6(32.19)&6(21.60)&4(21.30)&4(21.56)\\
  FAd-SGMRES-Sh&105(0.71)&60(0.46)&64(0.50)&63(0.46)&84(0.62)&70(0.49)&23(0.21)\\
  GMRES-DR-Sh&\multirow{2}{*}{$\dag$}&\multirow{2}{*}{$\dag$}&\multirow{2}{*}{$\dag$}&\multirow{2}{*}{$\dag$}&\multirow{2}{*}{$\dag$}&\multirow{2}{*}{$\dag$}&\multirow{2}{*}{$\dag$}\\
  $e=3$&&&&&&&\\
  GMRES-DR-Sh&\multirow{2}{*}{$\dag$}&\multirow{2}{*}{$\dag$}&\multirow{2}{*}{$\dag$}&\multirow{2}{*}{$\dag$}&\multirow{2}{*}{$\dag$}&\multirow{2}{*}{$\dag$}&\multirow{2}{*}{$\dag$}\\
  $e=6$&&&&&&&\\
  FAd-SGMRES-DR-Sh&\multirow{2}{*}{80(0.57)}&\multirow{2}{*}{54(0.42)}&\multirow{2}{*}{56(0.42)}&\multirow{2}{*}{57(0.42)}&\multirow{2}{*}{71(0.51)}&\multirow{2}{*}{52(0.39)}&\multirow{2}{*}{23(0.21)}\\
  $e=3$&&&&&&&\\
  FAd-SGMRES-DR-Sh&\multirow{2}{*}{74({\bf 0.50)}}&\multirow{2}{*}{51({\bf 0.35})}&\multirow{2}{*}{53({\bf 0.35})}&\multirow{2}{*}{56({\bf 0.37})}&\multirow{2}{*}{70({\bf 0.49})}&\multirow{2}{*}{48({\bf 0.31})}&\multirow{2}{*}{23({\bf 0.16})}\\
  $e=6$&&&&&&&\\
  \hline
\end{tabular}
\label{tab7}
\end{table}
 %\begin{tabular*}{\textwidth}{@{\extracolsep\fill}ccccc}
  %\hline
  %\multicolumn{5}{c}{$mv(cpu),~m=10,~\mu=0.9$}\\
  %\cline{2-5}
  %\multirow{2}{*}{Matrix}&\multirow{2}{*}{Ad-SGMRES-Sh}&\multirow{2}{*}{FAd-SGMRES-Sh}&\multicolumn{2}{c}{FAd-SGMRES-Dr-Sh}\\
  %\cline{4-5}
  %&&&$e=3$&$e=6$\\
  %\hline
  %$A_1$ &812(0.875) &105(0.625)      &80(0.500)      &74({\bf 0.485})  \\
  %$A_2$ &315({\bf 0.313})  &60(0.360) &54(0.344)       &51(0.328)  \\
  %$A_3$ &634(0.719)  &64(0.407)       &56(0.391)       &53({\bf 0.375})  \\
  %$A_4$ &384(0.406)  &63(0.359)       &57({\bf 0.344}) &56(0.36)        \\
  %$A_5$ &564(0.579)  &84(0.500)      &71({\bf 0.437})       &70(0.453)  \\
  %$A_6$ &2357(2.344) &70(0.406)       &52(0.312)       &48({\bf 0.297})  \\
  %$A_7$ &176(0.187)  &23({\bf 0.141}) &23(0.147)       &23(0.157)        \\
  %\hline
%\end{tabular*}

\begin{table}[!htpb]
\centering
\caption{Convergence behaviors of the Ad-SGMRES-Sh, FAd-SGMRES-Sh and FAd-SGMRES-DR-Sh with $n=49152$, ${\tt tol = 1e-6}$, $m=20$ and $\mu=0.9$}
\vspace{0.5mm}
 \begin{tabular}{cccccc}
  \hline
  \multicolumn{6}{c}{$mv(cpu),~m=20,~\mu=0.9$}\\
  \cline{2-6}
  \multirow{2}{*}{Matrix}&\multirow{2}{*}{Ad-SGMRES-Sh}&\multirow{2}{*}{FAd-SGMRES-Sh}&
  \multicolumn{3}{c}{FAd-SGMRES-DR-Sh}\\
  \cline{4-6}
  &&&$e=5$&$e=10$&$e=15$\\
  \hline
  $A_8    $&872(18.71)&105(12.11)     &95({\bf 11.27}) &94(11.70) &92(12.90) \\
  $A_9    $&584(12.90) &79(9.47)     &77({\bf 9.10})  &76(9.48)  &76(10.93) \\
  $A_{10} $&471(10.34) &72(8.46) &71(8.50)  &69({\bf 8.41})  &69(9.38)  \\
  $A_{11} $&431(9.61) &71({\bf 8.32})&72(8.66) &71(8.80)  &71(9.75)  \\
  $A_{12} $&659(15.11) &53(6.63)       &50(6.06)   &50({\bf 5.94})   &50(6.31)   \\
  $A_{13} $&1010(21.72)&54(6.27)       &51({\bf 5.94})   &52(6.18)   &51(6.61)   \\
  $A_{14} $&648(13.69) &54(6.13)       &49({\bf 5.63})   &49(5.83)   &49(6.19)   \\
  \hline
\end{tabular}
\label{tab8}
\end{table}

%\newpage
%\medskip
%From the tables, we can also observe that for most examples, the sum of the $mv$ %and $imv$
%of FAd-SGMRES-Sh and FAd-SGMRES-Dr-Sh are much larger than that of Ad-SGMRES-Sh, but the $cpu$ elapsed in FAd-SGMRES-Sh and FAd-SGMRES-Dr-Sh are much smaller than that of Ad-SGMRES-Sh. This is because that in the preconditioning process of FAd-SGMRES-Sh and FAd-SGMRES-Dr-Sh, %when using the matrix-vector product
%there is no solving process for the whole shifted linear systems, while for the matrix-vector product in the outer iteration, the solving process needs to run, thus there should be more time to elapse with more matrix-vector products in the outer iteration.

%\begin{flushleft}
%{\bf Table 8}
%\end{flushleft}
%\begin{tabular*}{\textwidth}{@{\extracolsep\fill}ccccc}
%\hline
%Matrix ID & Matrix name & Size   & Nonzeros & Problem domain               \\
%\hline
%1       & cnr-2000       & 325,557  & 3,216,152   & \\
%2       & email-EuAll     & 265,214  & 420,045    & Email network form a EU research institution            %\\
%3       & enron     & 69,244  & 276,143    & Enron email network                     \\
%4       & eu-2005       & 862,664    & 19,235    & \\
%5       & in-2004        & 1,382,908 & 16,917,053   & \\
%6       & indochina-2004    & 7,414,866  & 194,109,311    & \\
%7       & wb-cs-stanford       & 9,914  & 36,854   & Stanford CS web         \\
%  
%

{\bf Example 3.3}\quad As we know, preconditioning is the critical point that 
effects the convergence of iteration methods directly \cite{Saa2003}. However, 
different preconditioners will make different effects. In this example, some 
numerical results of FAd-SGMRES-Sh with different preconditioners are reported. 
We select ILU and IGMRES \cite{Saa2003}, and then denote the two algorithms by 
FAd-SGMRES-Sh(ILU) and FAd-SGMRES-Sh(IGMRES), respectively. At the same time, 
we also execute the flexible preconditioned GMRES with LU decomposition (FGMRES-Sh(LU)) 
\cite{GuZhoLin2007} for comparison. All the matrices used in the above two examples 
are considered in our experiments, and record the typical results in Table \ref{tab9}. 
Here $iter$ denotes the iteration number of Arnoldi process.

As seen from Table \ref{tab9}, FGMRES-Sh(LU) and FAd-SGMRES-Sh(ILU) are almost 
the same performance for each matrices. Especially for smaller size matrices, they 
are both performing better than FAd-SGMRES-Sh(IGMRES). However, for large-scale 
matrices, FAd-SGMRES-Sh(IGMRES) will be the best solver. This is because the inner 
loop of FGMRES-Sh(LU) becomes time-consuming to exactly solve a linear system with 
the coefficient matrix $A+\sigma_iI$ using the LU decomposition, and the saving of 
the LU decomposition is another big cost. For FAd-SGMRES-Sh(ILU), although there is 
no storage about the LU decomposition, but in each cycle, there needs to calculate 
the incomplete LU decomposition of $A$ and solving two sparse triangular linear 
systems, these are still both flaws. While for FAd-SGMRES-Sh(IGMRES), $10$ steps
of the inexact GMRES will not cost too much time. Consequently, for smaller size matrices, 
it is better to use FGMRES-Sh(LU) and FAd-SGMRES-Sh(ILU) to solve shifted systems, 
and it is best to use FAd-SGMRES-Sh(IGMRES) for solving some large-scale shifted systems.
\begin{table}[t]
\centering
\caption{Convergence behaviors of the FGMRES-Sh(LU), FAd-SGMRES-Sh(ILU) and 
FAd-SGMRES-Sh(IGMRES) with ${\tt tol = 1e-6}$, $m=20$, $\mu=0.9$, $\alpha=[0,0.4,2]$,
and $\sigma_1=0.5$, $\sigma_2=1$}
\vspace{0.5mm}
\begin{tabular}{cccc}
  \hline
  \multicolumn{4}{c}{$iter(cpu),~m=20,~\mu=0.9$}\\
  \cline{2-4}
  Matrix&FGMRES-Sh(LU)&FAd-SGMRES-Sh(ILU)&FAd-SGMRES-Sh(IGMRES)\\
  \hline
  bidiag1&8(0.33)&3({\bf 0.06})&42(0.09)\\
  sherman4&14({\bf 0.07})&17(0.16)&24(0.10)\\
  wang4&1181(199.68)&1241(2433.90)&117({\bf 1.88})\\
  young1c&11({\bf 0.15})&13(0.16)&299(0.41)\\
  young2c&10({\bf 0.03})&13(0.11)&265(0.31)\\
  \hline
\end{tabular}
\label{tab9}
\end{table}
%%%%%%%%%%%%%%%%%%%%%%%%%%%%%%%%%%%%%%%
\section{Conclusions}
\label{sec4}
\quad\
In the present paper, we established two iterative algorithms based on the Simpler 
GMRES for solving shifted linear systems simultaneously, namely FAd-SGMRES-Sh and 
FAd-SGMRES-DR-Sh. Moreover, these variants can be regarded as two improvements of 
Ad-SGMRES-Sh, which is recently proposed by Jing, Yuan and Huang in \cite{JinYuaHua2016}. 
The resultant algorithms converge in less matrix-vector products than the other related 
solvers (GMRES-Sh, Ad-SGMRES-Sh, FAd-GMRES-Sh, and GMRES-DR-Sh), especially for large 
problems. Furthermore, although the cost per iteration of FAd-SGMRES-Sh and FAd-SGMRES-DR-Sh
is higher, in our numerical experiences, the overall execution time is still lower. 
In addition, the FAd-SGMRES-DR-Sh performs better than FAd-SGMRES-Sh when the coefficient 
matrix of the seed system has many eigenvalues close to the origin as verified by 
numerical experiments. In conclusion, the proposed algorithms can be recommended as 
two efficient tools for solving shifted linear systems.

As an outlook for the future, the advanced development of preconditioning strategies 
(such as the polynomial preconditioning \cite{WuWanJin2012,Ahmad2017}, the nested 
iterative technique \cite{Baumann2015} and other preconditioning strategies \cite{
Simon2010,Soodx2016}) for solving shifted linear systems remains an meaningful topic 
of further research.
\section*{Acknowledgements}
{\em The authors would like to thank the anonymous referees and Editor-in-Chief Prof. 
Leszek Feliks Demkowicz for their constructive comments that helped to improve the 
quality of the paper.}

\def\refname
\end{document}